\documentclass[leqno,11pt]{amsart}

\usepackage{amssymb, amsmath}
\def\inte#1{
\displaystyle\mathop{#1\kern0pt}^\circ }

\def\supetage#1#2{
\sup_{\scriptstyle {#1}\atop\scriptstyle {#2}} }

\let\e=\varepsilon

\let\lam=\lambda

\def\cF{{\mathcal F}}

\def\virgp{\raise 2pt\hbox{,}}
\def\cdotpv{\raise 2pt\hbox{;}}

\def\eqdefa{\buildrel\hbox{{\rm \footnotesize def}}\over =}

\def\C{\mathop{\bf C\kern 0pt}\nolimits}
\def\DD{\mathop{\bf D\kern 0pt}\nolimits}
\def\K{\mathop{\bf K\kern 0pt}\nolimits}
\def\N{\mathop{\bf N\kern 0pt}\nolimits}
\def\Q{\mathop{\bf Q\kern 0pt}\nolimits}
\def\R{\mathop{\bf R\kern 0pt}\nolimits}
\def\SS{\mathop{\bf S\kern 0pt}\nolimits}
\def\ZZ{\mathop{\bf Z\kern 0pt}\nolimits}
\def\TT{\mathop{\bf T\kern 0pt}\nolimits}
\def\P{\mathop{\bf P\kern 0pt}\nolimits}

\def\dive{\mathop{\rm div}\nolimits}
\def\curl{\mathop{\rm curl}\nolimits}

\def\R{\mathop{\mathbb  R\kern 0pt}\nolimits}

\def\P{\mathop{\mathbb  P\kern 0pt}\nolimits}

\newcommand{\beq}{\begin{equation}}
\newcommand{\eeq}{\end{equation}}
\newcommand{\ben}{\begin{eqnarray}}
\newcommand{\een}{\end{eqnarray}}
\newcommand{\beno}{\begin{eqnarray*}}
\newcommand{\eeno}{\end{eqnarray*}}

\newtheorem{theo}{Theorem}

\newtheorem{rmk}{Remark}[section]

\newtheorem{prop}{Proposition}[section]

\begin{document}

\title[Remarks on the blow-up of a toy model for the Navier-Stokes equations]{Remarks on the blow-up  of solutions to  a   toy model for the Navier-Stokes equations}

\author[I. Gallagher]{Isabelle Gallagher}
\address[I. Gallagher]%
{ Institut de Math{\'e}matiques de Jussieu UMR 7586\\ Universit{\'e} Paris 7\\
175, rue du Chevaleret\\ 75013 Paris\\FRANCE }
\email{Isabelle.Gallagher@math.jussieu.fr}

\author[M. Paicu]{Marius Paicu}
\address[M. Paicu]%
{ D\'epartement de Math\'ematiques\\ Universit{\'e} Paris 11\\
B\^{a}timent 425\\ 91405 Orsay Cedex\\FRANCE }
\email{marius.paicu@math.u-psud.fr }


\begin{abstract}
In~\cite{ms}, S. Montgomery-Smith provides a one dimensional model for the three dimensional, incompressible Navier-Stokes 
equations, for which he proves the 
blow up of solutions associated to a class of large initial data, while the same global existence results as for the Navier-Stokes equations hold for small data. In this note the model is adapted to the case of two and three
space dimensions, with the additional feature that the divergence free condition is preserved. It is checked that the family of initial data  constructed in~\cite{cg}, which is arbitrarily large but yet generates a global solution to the Navier-Stokes equations in three space dimensions, actually causes blow up for the toy model --- meaning that the precise structure of the nonlinear term is crucial to
understand the   dynamics of large solutions to the  Navier-Stokes equations. 
\end{abstract}
  

\maketitle


\section{Introduction}
Consider the Navier-Stokes equations in~$\R^{d}$, for~$d = 2$ or~3,
 \[
{\rm (NS)}\ \left\{
\begin{array}{c}
\partial_{t} u -\Delta u+u\cdot\nabla u =-\nabla p\\
\mbox{div}\:   u =0 \\
u_{|t = 0} = u_{0},
\end{array}
\right.
\]
where~$u = (u^1,\dots,u^d)$ is the velocity of an incompressible, viscous, homogeneous fluid evolving in~$\R^d$, and~$p$ is its pressure. Note that the divergence free condition  allows to recover~$p$ from~$u$ through the formula
$$
-\Delta p = \mbox{div} \: (u\cdot\nabla u).
$$
A formally equivalent formulation for~(NS) can be obtained
by applying the  projector onto divergence free vector fields~$\displaystyle
\P \eqdefa \mbox{Id} - \nabla
\Delta^{-1} \mbox{div}  
$ to~(NS):
 \[
  \left\{
\begin{array}{c}
\partial_{t} u -\Delta u+\P  (u\cdot\nabla u) =0\\
u_{|t = 0} = u_{0} = \P u_0.
\end{array}
\right.
\]

This system has three important features: 

(E) (the {\it energy} inequality): the~$L^{2} (\R^d)$ norm of~$u$ is formally bounded for all times by that of the initial data;

(I) (the {\it incompressibility} condition): the solution  satisfies for all times the constraint~$\mbox{div}\:   u =0 $; 

(S) (the {\it scaling} conservation): if~$u$ is a solution associated with the data~$u_0$, then for any positive~$\lambda$, the rescaled~$u_{\lam}(t,x) \eqdefa \lam u(\lam^{2}t, \lam x)$
 is a solution associated with~$u_{0,\lam}(x) \eqdefa \lam u_{0}( \lam x)$.

Of course the two first properties are related, as (I) is the ingredient enabling one to obtain~(E), due to the special structure of the nonlinear term.

\medskip

Taking~(E) into account, one can prove   the existence of global, possibly non unique, finite energy solutions (see the fundamental work
 of J. Leray~\cite{leray}). On the other hand the use of~(S) and a fixed point argument enables one to prove the existence of a unique, 
 global solution if the initial data is small  
in scale-invariant spaces (we will call ``scale-invariant space'' any Banach space~$X$ satisfying~$\|\lam f(\lam \cdot)
\|_{X} = \|f\|_{X}$ for all~$\lam>0$): for instance the homogenenous Sobolev space~$\dot H^{\frac d2 - 1}$, 
Besov spaces~$\dot B^{-1 + \frac dp}_{p,\infty}$ for~$p<\infty$ or the  space~$BMO^{-1}$. We recall that
 $$
 \|f\|_{\dot B^{s}_{p,q}} \eqdefa \left\| t^{-\frac s2} \|e^{t\Delta}f\|_{L^p(\R^d)} \right\|_{L^q(\R^+;\frac{dt}t)},
 $$
 and
$$
\|f\|_{BMO^{-1}} \eqdefa \sup_{t > 0}
\biggl( t^{\frac 1 2}\|e^{t\Delta} f\|_{L^\infty}
+\supetage{x\in \R^d}{R>0} R^{-\frac d 2}
\Big(\int_{P(x,R)}| e^{t\Delta} f (t,y)|^2 dy\Bigr)^{\frac 1 2}\biggr)  ,
$$
where~$P(x,R)=[0,R^2]\times B(x,R)$
and~$B(x,R)$ denotes the ball  of~$\R^d$ of center~$x$ and radius~$R$.
We refer 
respectively to~\cite{fk},\cite{cmp} and~\cite{kt} for proofs of the wellposedness of~(NS) for small data in those spaces. When~$d=2$, the smallness condition may be removed: that has been  known since the work of J. Leray (\cite{Leray2D}) in the energy space~$L^2$ (which is scale invariant in two space dimensions), and was proved in~\cite{gp},\cite{g} for larger spaces, provided they are completions of the Schwartz class for the corresponding norm (Besov   or~$BMO^{-1}$ norms).

It is well known and rather easy to see that the largest scale invariant Banach space embedded in the space of tempered distributions is~$\dot B^{-1  }_{\infty,\infty}$. In three or more space dimensions, it is not known that global solutions exist for smooth data, arbitrarily large in~$\dot B^{-1  }_{\infty,\infty}$. We will not review here all the progress made in that direction in the past years, but merely recall a few of the main recent achievements concerning the possibility of blow up of large solutions. Recently, D. Li and Ya. Sinai were able in~\cite{lisinai} to prove the blow up in finite time of solutions to the  Navier-Stokes equations for complex initial data. We note that, as for the system that we construct in the present paper, the complex Navier-Stokes system does not satisfy  any energy inequality. Before that, some numerical evidence was suggested to support the idea of finite time blow up of (NS) (see for instance~\cite{no} or~\cite{grundymclaughlin}).
On the other hand in~\cite{cg}  a class of 
large initial data was constructed, giving rise to a global, unique solution; this family will be presented below.  Another type of  example was provided in~\cite{cg2}. It should be noted that in both those examples, the special structure of the equation is crucial to obtain the global wellposedness. In~\cite{ms}, S. Montgomery-Smith suggested a   model for~(NS), with the same scale invariance and for which
the same global wellposedness results hold for small data. The interesting feature of   the model is that it is possible (see~\cite{ms}) to prove the blow up in finite
time of
some solutions. The model is the following:
$$
{\rm (TNS}_{1}{\rm )} \quad \quad  \partial_{t} u -\Delta u = \sqrt{-\Delta} \: (u^{2})    \quad \mbox{in} \: \R^{+} \times \R.
$$
The main ingredient  of the proof of the existence of blowing-up solutions consists in noticing that if the initial data  has a positive Fourier transform, then that positivity is preserved for the solution at all further times. One can then use the Duhamel formulation of the solution and deduce a lower bound for the Fourier transform that blows up in finite time. We will not  write more details here as we will be reproducing that computation in Section~\ref{2D}.

\medskip

In this paper we adapt the construction of~\cite{ms} to higher space dimensions. In order to have a proper   model
in higher dimensions it is important to preserve as many features of the Navier-Stokes equations as possible. Here we will seek to preserve scaling~(S) as well as  the divergence free condition~(I) (as we will see, condition~(E) cannot be preserved in our model). This amounts to transforming the nonlinear
term proposed in~\cite{ms} (see Equation~(TNS$_{1}$)  above) in such a way as to preserve both the positivity conservation property  in Fourier space {\it and} the incompressibility condition. This
is in fact a technicality which may be handled by explicit computations in Fourier space; actually the more interesting aspect of the result we obtain is that the initial data constructed in~\cite{cg} to show the possibility of global solutions associated with arbitrarily large initial data actually generates a blow-up solution for~(TNS$_{3}$). This, joint to the fact that we are also able to 
obtain blowing-up solutions in the two dimensional case, indicates that proving a global existence result for arbitrarily large data for~(NS) requires using the energy estimate, or the specific structure of the nonlinear term --   two properties which are discarded in our model.

Let
us state  the result proved in this paper.

\begin{theo}\label{ms2D3D}
Let the dimension~$d $ be equal to 2 or 3. There is a bilinear operator~$Q$, which is a~$d$-dimensional  matrix of Fourier multipliers of order one, such 
as the equation
$$
{\rm (TNS}_{d}{\rm )}  \quad  \left\{
\begin{array}{c}
\partial_{t} u -\Delta u = Q(u,u)   \quad \mbox{in} \: \R^{+} \times \R^{d}\\
\mbox{div}\:   u =0 \\
u_{|t = 0} = u_{0}
\end{array}
\right.
$$
satisfies  properties (I) and~(S), and such that   there is a global, unique solution if the data
is small enough in~$BMO^{-1}$. Moreover there is a family of smooth initial data~$u_{0}$, which may be chosen arbitrarily large
in~$ \dot B^{-1  }_{\infty,\infty}$, such that the associate solution of~(TNS$_{d}$) blows up in all Besov norms, whereas the 
associate solution of~(NS) exists globally in time.
\end{theo}
The proof of the theorem is given in the sections below. In Section~\ref{2D}
we deal with the two dimensional case, while the three dimensional case is treated in  Section~\ref{3D}: in both cases
we  present an alternative to the   bilinear term of (NS), which preserves scaling  and the divergence free property, while giving rise to  solutions blowing up in finite time, for some classes of initial data. The fact that some of those initial data in fact generate a global solution for  the three dimensional Navier-Stokes equations is addressed in   Section~\ref{exampleC-1}.
\begin{rmk}
We note that the method of the proof allows to construct  blowing up   solutions for the hyper-viscous case,   meaning for equations of the form
$$
 \quad  \left\{
\begin{array}{c}
\partial_{t} u -\Delta^{\alpha} u = Q(u,u)   \quad \mbox{in} \: \R^{+} \times \R^{d}\\
\mbox{div}\:   u =0 \\
u_{|t = 0} = u_{0}
\end{array}
\right.
$$
where $\alpha\geq 1$ and $\Delta^{\alpha}f=\cF^{-1}(|\xi|^{2\alpha}\hat f(\xi))$. Indeed, the only important feature  in order to construct  blowing-up  solutions by the method of \cite{ms} is that the system written in the Fourier variable, preserves the positivity of the symbol, and so, the positivity of $\hat u^j(t,\xi)$ if $\hat u_0^j(\xi)>0$, for any~$j \in \{1, \dots, d\}$. 
\end{rmk}

\begin{rmk} In the two dimensional case, it    might  seem  more natural to work on the vorticity formulation of the equation: in 2D it is well known that the vorticity satisfies a transport-diffusion equation, which provides easily the existence of global solutions for any sufficiently smooth initial data. An example where the vorticity equation is modified (rather than~(NS)) is provided at the end of Section~\ref{2D} below.
\end{rmk}


\section{Proof of the theorem in the two-dimensional case}\label{2D}
 In this section we shall construct the quadratic form~$Q$, as given in the statement of Theorem~\ref{ms2D3D}, which allows to construct  blowing up  solutions for the (TNS$_2$) system. 

Let us consider a system of the following form:
$$
 \quad  \left\{
\begin{array}{c}
\partial_t u - \Delta u= {\mathcal Q}(u , u )-\nabla p  \\
\dive u=0.
\end{array}
\right.$$
 Taking the Leray projection of this equation, we obtain 
$$\partial_t u-\Delta u=\P {\mathcal Q} (u ,u ).$$
We wish to follow the idea of the proof of~\cite{ms}, thus to construct~$Q = \P {\mathcal Q} $ as a matrix of   Fourier multipliers of order~1, such that the product~$\widehat{\P   {\mathcal Q} }$ preserves the positivity of the Fourier transform.
We define~$Q (u,u)$ as the vector whose $j$-component is, for~$ j \in \{1,2\},$
\beq \label{defQ}
  \left(Q (u , u ) \right)^j = \sum_i q_{i,j} (D) (u^i u^j),
\eeq
and we impose that~$q_{i,j} (D)$
 are Fourier multipliers of order~1.
For example, let us simply choose 
$$
\widehat   {\mathcal Q} (\xi) = |\xi|
 {\mathbf 1}_{ \xi_{1} \xi_2<0}   \left(
\begin{array}{cc}
1&1\\
1&1 \\
\end{array}
\right).
$$
Recalling that 
$$
\widehat \P(\xi) = 
\left(
\begin{array}{cc}
1-\frac{\xi_1^2}{|\xi|^2}&-\frac{\xi_1\xi_2}{|\xi|^2}\\
-\frac{\xi_1\xi_2}{|\xi|^2}&1-\frac{\xi_2^2}{|\xi|^2} \\
\end{array}
\right),
$$
we easily obtain 
$$
\widehat{\P   {\mathcal Q} }(\xi) = 
 {\mathbf 1}_{ \xi_{1} \xi_2<0}\frac{1}{|\xi| }\left(
\begin{array}{cc}
\xi_2^2-\xi_1\xi_2&\xi_2^2-\xi_1\xi_2\\
\xi_1^2-\xi_1\xi_2&\xi_1^2-\xi_2\xi_1\\
\end{array}
\right),
$$
so all the elements of this matrix are positive. 

\medskip

The Duhamel formulation of~(TNS$_2$) reads
$$
\widehat u^j(t,\xi)=e^{-t|\xi|^2}\widehat u_0^j(\xi)+ \sum_i \int_0^t e^{-(t-s)|\xi|^2}q_{i,j}(\xi)(\widehat u^i(s)\ast \widehat u^j (s)) \: ds,
$$
where we have denoted by~$q_{i,j}(\xi)$ the matrix elements of~$\widehat{\P  {\mathcal Q}}(\xi)$.

 It  is not difficult to see that all the usual results on the Cauchy problem for the Navier-Stokes equations hold for this system (namely results of~\cite{fk},\cite{cmp} and~\cite{kt} as recalled in the introduction).
 Moreover it 
is clear that if the Fourier transform of~$\widehat u_0$ is positive, then that positivity property holds for all times. 
 
  Now let us construct a data generating a   solution blowing up in finite time. We will be following closely  the argument of~\cite{ms}, and we refer to that article for all the computational details. We start  by choosing the  initial data~$u_0=(u_0^1, u_0^2)$ such that~$\widehat u_0^1 \geq 0$, and the support of $\widehat u_0^1$ lies in the second and fourth sector of the complex plane, that is the zone where $\xi_1\xi_2<0$; we also suppose
  this spectrum is symmetric  with  respect to zero (and to fix notation, that  the support of~$\widehat u_0^1 $ intersects the set~$|\xi_j| \geq 1/2$, for~$j \in \{1,2\}$).  Taking into account the   divergence free  condition  which states that $\displaystyle \widehat u^2(\xi)=-\frac{\xi_1 \xi_2 }{\xi_2^2 }\widehat u^1(\xi) $, we deduce that $\widehat u_2$ is supported in the same region as $\widehat u_1 $ and is also nonnegative.
  
   Let us denote by~$A$ the~$L^1$ norm of~$u_0$ (which will be assumed to be large enough at the end), and let us write~$u_0 = Aw_0$.   
  The idea, as in~\cite{ms}, is to prove   that for any~$k \in \N$ and~$j \in \{1,2\}$, 
\beq \label{induction}
  \widehat u^j (t,\xi) \geq A^{2^{k}} e^{-2^k t} 2^{k-4(2^k-1)} {\mathbf 1}_{t \geq t_k}  \widehat w_0^{k,j} (\xi)
  \eeq
where we have written, $  w_0^{k,j} = (w_0^{0,j})^{2^k} $
and~$ \widehat w_0^0$ is the restriction of~$ \widehat w_0 \:  {\mathbf 1}_{|\xi_j| \geq 1/2}$ to the second sector of the plane. Finally the time~$t_k$
is chosen so that~$t_0 = 0$ and~$t_k - t_{k-1} \geq 2^{-2k} \log 2$. Notice that~$\displaystyle \lim_{k \rightarrow \infty} t_k = \log 2^{1/3}. $  The result~(\ref{induction}) is proved by induction. Suppose that~(\ref{induction}) is true for~$k-1$ (it is clearly true for~$k = 0$). Due to the positivity of~$\widehat u_0$, we can write
\begin{eqnarray*}
\widehat u^j(t,\xi) & \geq&  \sum_i  \int_0^t e^{-(t-s)|\xi|^2}q_{i,j}(\xi)(\widehat u^i(s,\xi)\ast \widehat u^j (s,\xi)) \: ds \\
 & \geq&  \int_0^t e^{-(t-s)|\xi|^2}q_{j,j}(\xi)(\widehat u^j(s,\xi)\ast \widehat u^j (s,\xi)) \: ds 
 \end{eqnarray*}
and using the induction assumption, along with the support restriction  of~$w_0^{k-1,j}$, we find that
\begin{eqnarray*}
 \widehat u^j(t,\xi) 
 & \geq&  \int_0^t e^{-(t-s)|\xi|^2}q_{j,j}(\xi)  (A^{2^k-1} \alpha_{k-1}(s))^2   \: ds \:   \widehat w_0^{k-1,j} \ast \widehat w_0^{k-1,j}(\xi)\\
& \geq&  \int_0^t e^{-(t-s)2^{2k}}q_{j,j}(\xi)  (A^{2^k-1} \alpha_{k-1}(s))^2   \: ds \: \widehat w_0^{k-1,j} \ast \widehat w_0^{k-1,j}(\xi),
 \end{eqnarray*}
 where~$\alpha_k(t) =2^{k-4(2^k-1)} {\mathbf 1}_{t \geq t_k}  $. But $ \widehat w_0^{k-1,j} \ast \widehat w_0^{k-1,j} =  \widehat w_0^{k ,j}$, and on the support of~$\widehat w_0^{k ,j}$
we have~$q_{j,j}(\xi) \geq C 2^k$. The induction then follows exactly as in~\cite{ms}. 

Once~(\ref{induction}) is obtained, the blow up of all~$\dot B^{s}_{\infty,\infty}$ norms follows directly, noticing that~$u^j(t_\infty)$ can be bounded from below in~$\dot B^{s}_{\infty,\infty}$ by~$C (Ae^{-t_\infty} 2^{-4})^{2^k} 2^{(s+1)k}$, which goes to infinity with~$k$ as soon as~$Ae^{-t_\infty} 2^{-4} >1$. That lower bound is simply due to the fact that (calling~$\Delta_k$ the usual Littlewood-Paley truncation operator  entering in the definition of Besov norms)
$$
\|u (t_\infty)\|_{\dot B^{s}_{\infty,\infty}} = \sup_k 2^{ks} \|\Delta_k u (t_\infty) \|_{L^\infty} \geq  \sup_k  2^{ks} |\Delta_k u (t_\infty,0)| =  \sup_k  2^{ks}
\|\widehat {\Delta_k u} (t_\infty) \|_{L^1}
$$
since~$\widehat {\Delta_k u}(t_\infty)$ is nonnegative.

\begin{rmk} \label{whyitworks}
One can notice that as soon as the  matrix~$Q$ has been defined, the computation turns out to be  identical to the case studied in~\cite{ms}.   In particular the important fact is that~$\widehat u_0$ is nonnegative (and that its support intersects, say, the set~$|\xi_j| \geq 1/2$).
\end{rmk}

\begin{rmk}  As explained in the introduction, it seems natural to try to improve the previous example by perturbing the vorticity equation, since that equation is   special in two space dimensions.  
Let us therefore consider   the vorticity $\omega=\partial_1 u^2-\partial_2 u^1$. As is well known, the two dimensional Navier-Stokes equations can simply be written as a transport-diffusion equation on~$\omega$:  
$$\partial_t\omega+u \cdot \nabla\omega-\Delta \omega=0,$$
which can also be written, since~$u$ is divergence free,
$$
\partial_t\omega + \partial_1 (u^1\omega) + \partial_2(u^2\omega)- \Delta\omega=0.
$$

Changing the place of the derivatives, and noticing that a derivative of~$u$ has the same scaling as~$\omega$, a model equation for the vorticity equation is simply
$$
\partial_t\omega +\omega^2 - \Delta\omega=0.
$$
This simplified model is a semilinear heat equation for which the blow-up of the solution is  well known (see \cite{friedman}, \cite{fujita}). It is also easy to see that the argument of~\cite{ms} is true for this system, which therefore blows up in finite time for large enough initial data with negative Fourier transform. One can note that the equation on $u$ becomes
$$
\partial_t u+\nabla^\perp\Delta^{-1}\big((\curl u)^2\big)-\Delta u=-\nabla p \quad,\quad \dive u=0,
$$
which blows up but does not preserve the sign of the Fourier transform.

\end{rmk}


\section{Proof of the theorem in the three-dimensional case}\label{3D}
The three-dimensional situation follows the lines of the two-dimensional case studied above, though it is slightly more technical.  The main step, as in the previous section, consists in finding a three-dimensional matrix~${\mathcal Q}$ such that the Fourier transform of the product~$\P{\mathcal Q}$ has positive coefficients (we recall that~$\P$ denotes the~$L^2$ projection onto divergence free vector fields).
Let us define, similarly to   the previous section,  the matrix
$$
\widehat {\mathcal Q} (\xi) = |\xi|   {\mathbf 1}_{ \xi \in {\mathcal E}} \left(
\begin{array}{ccc}
1&1&1 \\
1&1&1 \\
1&1&1 
\end{array}
\right), $$
where~$\displaystyle
 {\mathcal E} \eqdefa \left\{\xi \in \R^3, \: \xi_{1}\xi_{2} <0, \: \xi_{1}\xi_{3} <0, \: 
 | \xi_{2}| < \min(|\xi_1 |,
 |\xi_3 |)\right\}.$
We compute easily that
$$ \widehat {\P {\mathcal Q} }(\xi) = 
  {\mathbf 1}_{ \xi\in {\mathcal E}}   |\xi|^{-1} \left(
\begin{array}{ccc}
\xi_2^2 + \xi_3^2 - \xi_1 \xi_2 - \xi_1 \xi_3 &\xi_2^2 + \xi_3^2 - \xi_1 \xi_2 - \xi_1 \xi_3&\xi_2^2 + \xi_3^2 - \xi_1 \xi_2 - \xi_1 \xi_3 \\
\xi_1^2 + \xi_3^2 - \xi_1 \xi_2 - \xi_2 \xi_3&\xi_1^2 + \xi_3^2 - \xi_1 \xi_2 - \xi_2 \xi_3&\xi_1^2 + \xi_3^2 - \xi_1 \xi_2 - \xi_2 \xi_3 \\
\xi_1^2 + \xi_2^2 - \xi_1 \xi_3 - \xi_2 \xi_3&\xi_1^2 + \xi_2^2 - \xi_1 \xi_3 - \xi_2 \xi_3&\xi_1^2 + \xi_2^2 - \xi_1 \xi_3 - \xi_2 \xi_3\end{array}
\right).
$$
Let us consider the sign of the matrix elements of~$\widehat {\P {\mathcal Q} }(\xi) $.
The first line of the above matrix is clearly made of positive scalars, due to the sign condition imposed on the components of~$\xi$. 
The components of the second line may be written
 $$
 \xi_1^2 + \xi_3^2 - \xi_1 \xi_2 - \xi_2 \xi_3 = \xi_1^2 - \xi_1 \xi_2 + \xi_3 (\xi_3 - \xi_2),
 $$
 which is also positive since either~$ \xi_2 $ and~$\xi_3$ are both positive, in which case~$\xi_3 > \xi_{2}$, or they are both negative in which case~$\xi_3 < \xi_{2}$.
 Similarly one has
 $$
 \xi_1^2 + \xi_2^2 - \xi_1 \xi_3 - \xi_2 \xi_3 =  \xi_1^2 + \xi_2^2- \xi_3 (\xi_1 + \xi_2) ,
 $$
 and either~$\xi_1>0$, $\xi_2<0$, $\xi_3<0$ and~$\xi_1+\xi_2>0$, or~$\xi_1<0$, $\xi_2>0$, $\xi_3>0$ and~$\xi_1+\xi_2<0$. So the third line is also made of positive real numbers.

  Now that it has been checked that all coefficients are positive, we just have to follow again the proof of the two dimensional case to obtain the expected result, showing the blow up of solutions to
  $$
  \partial_t u - \Delta u =  Q(u,u), \quad u_{| t = 0} = u_0,
  $$
  where~$ Q(u,u) = \P {\mathcal Q}(u,u)$ is the vector defined as in~(\ref{defQ}).
   We will not write all the details, which are   identical to the two-dimensional case, but simply give the form of the initial data, which is summarized in the next proposition.
   \begin{prop}\label{data3D}
    Let~$u_0$ be a smooth, divergence free vector field such that the components of~$\widehat u_0$ are even, nonnegative functions, such that the support of~$\widehat u_0$ intersects the set~$|\xi_j| \geq 1/2$, for~$j \in \{1,2,3\}$. Then the unique solution to~(TNS$_3$) associated with~$u_0$ blows up in finite time, in all Besov spaces.
    \end{prop}
   We will not detail the proof of that proposition, as it is identical to the two dimensional case (thus in fact to~\cite{ms}).

   Of course one must check that such initial data exists. The simplest  way to construct such an initial data is simply    to suppose it only has two nonvanishing components, say~$u_0^1$ and~$u_0^2$, and that the Fourier transform of~$u_0^1$ is supported in~${\mathbf 1}_{\xi_1 \xi_2 <0}$ while intersecting the set~$|\xi_j| \geq 1/2$. The divergence free condition ensures that the same properties hold for~$u_0^2$ (and~$u_0^3$ is assumed to vanish identically). An explicit example is provided in the next section.
   
     That ends the proof of the ``blowing up'' part of the theorem.  
  
  \begin{rmk}
  Notice that in that example, the energy inequality~(E) cannot be satisfied, as it would require that~$(Q(u,u) | u)_{L^2} \geq 0$, which cannot hold in our situation  if the Fourier transform of~$\widehat u$ is nonnegative.
  \end{rmk}
  
 \section{Examples of arbitrarily large initial data providing a blowing up solution  to~(TNS$_3$) and a global solution  to~(NS)}\label{exampleC-1}
 In this short section, we check that the initial data provided in \cite{cg} and which allows to obtain large, global solutions for Navier-Stokes equations, gives rise to a solution blowing up  in finite time for the modified three dimensional Navier-Stokes equation constructed in the previous section.
 
 More precisely we have the following result.
 \begin{prop}\label{NSvsTNS}
 Let~$\phi $ be a function in~${\mathcal S}(\R^3)$, such that~$\widehat\phi\geq 0$, and such that~$\widehat\phi$ is even and has its support  in the region~$ {\mathbf 1}_{ \xi_{1}\xi_{2} <0}  $, while intersecting the set~$|\xi_j| \geq 1/2$, for~$j \in \{1,2,3\}$.  Let~$\e$ and~$\alpha$ be given in~$ ]0,1[$, and consider the family of initial data
 $$u_{0,\varepsilon}(x) = (\partial_{2}\varphi_{\varepsilon}(x),
-\partial_{1}\varphi_{\varepsilon}(x), 0)
$$
where 
 $$
\varphi_{\varepsilon}(x) = \frac{({-\log \varepsilon})^{\frac15}}
{\e^{1-\alpha}} \cos  \left(
\frac{x_{3}}{\e}\right)
(\partial_1\phi)\Bigl(x_{1}, \frac{x_{2}}{\varepsilon^{\alpha}}, x_{3}\Bigr).
$$
Then for~$\e>0$ small enough, the unique solution of~(NS) associated with~$u_{0,\varepsilon}$ is smooth and global in time, whereas the  unique solution of~(TNS$_3$) associated with~$u_{0,\varepsilon}$ blows up in finite time, in all Besov norms.
\end{prop}
\begin{rmk} It is proved in~\cite{cg} that such initial data has a large~$\dot B^{-1}_{\infty,\infty}$ norm, in the sense that there is a constant~$C$ such that
$$
C^{-1} ({-\log \varepsilon})^{\frac15} 
 \leq \|u_{0,\varepsilon}\|_{\dot B^{-1}_{\infty,\infty}}
  \leq C ({-\log \varepsilon})^{\frac15}.
$$
\end{rmk}
 To prove  Proposition~\ref{NSvsTNS}, we notice that the initial data given in the proposition
is a particular case of the family of initial data presented in~\cite{cg}, Theorem~2, which generates a unique, global solution as soon as~$\e$ is small enough (in~\cite{cg} there is no restriction on the support of the Fourier transform and~$\partial_1 \phi$ is simply~$\phi$). So we just have to check that the initial data fits with the requirements of Section~\ref{3D} above, and more precisely that it satisfies the assumptions of Proposition~\ref{data3D}.
 Notice that 
$$
\widehat   \varphi_\varepsilon   (\xi)=
 \frac{({-\log \varepsilon})^{\frac15}}
{2\e^{1-2\alpha} } 
\left(i\xi_1\widehat\phi(\xi_1, \varepsilon  ^\alpha \xi_2,\xi_3+\frac 1\varepsilon  )+i\xi_1\widehat\phi(\xi_1, \varepsilon  ^\alpha \xi_2,\xi_3-\frac 1\varepsilon  )\right)
$$
We need to check that $\widehat u^i_{0,\varepsilon  }\geq 0$, for $i \in \{1,2,3\}$, and that the Fourier support  intersects the set~$|\xi_j| \geq 1/2$.
 We have
$$\widehat u_{0,\varepsilon  }(\xi)=
 \frac{({-\log \varepsilon})^{\frac15}}
{2\e^{1-2\alpha} } 
\left(-\xi_1\xi_2\widehat\phi(\xi_1,\varepsilon  ^\alpha\xi_2,\xi_3\pm\frac 1\varepsilon  ), \xi_1^2\widehat\phi(\xi_1,\varepsilon  ^\alpha\xi_2,\xi_3\pm\frac 1\varepsilon  ),0\right),$$
and we have clearly  the desired properties.

This ends the proof of the proposition, and of the theorem.


\end{document}